\documentclass{emsprocart}

\contact[jean.bertoin@math.uzh.ch]{Jean Bertoin, Institut f\"ur Mathematik, 
Universit\"at Z\"urich, 
Winterthurerstrasse 190, 
CH-8057 Z\"urich}

\def \\ { \cr }
\def\R{\mathbb{R}}
\def \1{1 \mkern -6mu 1} 
\def\N{\mathbb{N}}

\def\R{\mathbb{R}}

  \def \dd{{\rm d}}
\def \ee{{\rm e}}

\newtheorem{theorem}{Theorem}[section]

\theoremstyle{definition}

\title[Coagulation with limited aggregations]{Coagulation with limited aggregations}

\author[Jean Bertoin]{Jean Bertoin
}

\begin{document}

\begin{abstract}
Smoluchowski's coagulation equations can be used as elementary mathematical models for the formation of polymers. We review here some recent contributions on a variation of this model in which the number of aggregations  for each atom is {\it a priori} limited. 
Macroscopic results in the deterministic  setting can be explained  at the microscopic level by considering  a version of stochastic coalescence with limited aggregations, which can be related to the so-called random configuration model of random graph theory. 
\end{abstract}



\maketitle

\section{Introduction} 
In 1916, Smoluchowski \cite{Sm} introduced an infinite system of differential equations which are meant to describe the evolution of the concentrations of particles in a medium such that  pairs of particles merge as time passes.
Since then, this area has been intensively studied by physicists, chemists, and mathematicians. See for instance the monographs by Drake \cite{Drake} and by Dubovski \cite{Dubo}, and the survey by Aldous \cite{Aldous}; of course  not to mention the numerous references therein and further recent works.

In the original model of Smoluchowski, particles are entirely determined by their sizes. Here, 
we will be interested in the case when they have also arms (or stubs) which serve to perform aggregations.
In other words, each particle has initially  a certain number of potential links
which are consumed when this particle takes part to a coagulation event. The effect is that the total number of aggregations involving a given particle is limited by its initial number of arms,
and we shall see that this has a significant impact on concentrations.

The purpose of this text  is to survey some recent contributions to this topic, including explicit resolution and  discussion  of phenomena such as gelation and self-organized criticality. Formulas for the limiting concentrations when time goes to infinity hint at connexions with certain branching processes. These are explained by the study of the microscopic model. The latter  is a close relative to the so-called random configuration model, a random multigraph in which the degree sequence is pre-described.

The rest of this paper is organized as follows. Section 2 contains some background on Smoluchowski's equations, with a focus on the multiplicative kernel and the phenomenon of gelation. Section 3 is devoted to the system of coagulation equations with limited aggregations; specifically we shall discuss its explicit solutions, another phenomenon of gelation, and the existence of non-degenerate limiting concentrations. Section 4 deals with the microscopic model of stochastic coalescent with limited aggregations and its relation with random configurations  in the subcritical case when gelation does not occur. Finally we shall briefly present some recent results of Merle and Normand  in the super-critical case,  in particular an interesting phenomenon of self-organized criticality which arises in this setting. 

\section{Smoluchowski's coagulation equation in a nutshell}

\subsection{General setting}
Consider a medium in which particles coalesce pairwise as time passes. Typically,  each particle has a mass (or size), say $m>0$, and when two particles with masses $m$ and $m'$ coalesce, they are replaced by a single particle with mass $m+m'$. 
Smoluchoswki proposed a simple mean-field model for the evolution of the concentration $c_t(m)$ of particles of mass $m$ at time $t\geq0$. Assume for simplicity that masses only take integer values. The dynamics are governed by a symmetric kernel
$\kappa: \N^*\times \N^*\to [0,\infty)$ which specifies the rate at which two particles with masses  $m$ and $m'$ coalesce. Smoluchowski's coagulation equation reads for every $m\in\N^*$
\begin{equation}\label{Sm}
{\partial_t}c_t(m) =\frac{1}{2}\sum_{m'=1}^{m-1} 
c_t(m')c_t(m-m')\kappa(m',m-m')-c_t(m)\sum_{m'=1}^{\infty}c_t(m')\kappa(m,m').
\end{equation}
In the right hand side, the first term accounts for the creation of particles of mass $m$ as result of the coagulation of a particle of mass $m'<m$ and a particle of mass $m-m'$ (and the factor $1/2$ is due to an obvious symmetry), and the second term  for the disappearance of particles of mass $m$ after a coagulation with another particle. 

It is often convenient to rewrite \eqref{Sm} in the equivalent form
\begin{equation}\label{Sm'}
\partial_t\langle c_t,f\rangle = \frac{1}{2}\sum_{m=1}^{\infty} \sum_{m'=1}^{\infty}(f(m+m')-f(m)-f(m'))\kappa(m,m')c_t(m)c_t(m')\,,
\end{equation}
where $f:\N^*\to \R$ has finite support and $\langle c_t,f\rangle=\sum_{m=1}^{\infty}f(m)c_t(m)$. 
The informal description of the model might suggest that the average mass of particles in the medium
$$\langle c_t, {\rm Id}\rangle= \sum_{m=1}^{\infty} mc_t(m)
$$
should be conserved as time passes, and indeed if we could take $f={\rm Id}$ in \eqref{Sm'}, then we would get $\partial_t\langle c_t, {\rm Id}\rangle=0$. However this is not always true; more precisely  the {\it gelation time} 
$$T_{\rm gel}=\inf\{t> 0: \langle c_t, {\rm Id}\rangle\neq \langle c_0, {\rm Id\rangle }\}$$
is known to be finite for a large class of coagulation kernels (see e.g. \cite{EMP, ELMP}).

Informally, the phenomenon of gelation should be interpreted as the formation of giant particles which are not taken into account in the average mass of particles in the medium. 
The prototype of kernels for which this  occurs is the multiplicative one; it is  discussed next.

\subsection{The multiplicative kernel and gelation phenomenon}
We suppose here that $\kappa(m,m')=m\cdot m'$, and further, 
with little loss of generality,  that the average mass of particles in the media at the initial time $t=0$  is $\langle c_0, {\rm Id}\rangle=1$. Then for $t<T_{\rm gel}$,  Smoluchowski equation \eqref{Sm} takes the simpler form
\begin{equation}\label{Smolsimple}
{\partial_t} c_t(m)=\frac{1}{2}\sum_{m'=1}^{m-1} 
c_t(m')c_t(m-m')m' (m-m')-mc_t(m).
\end{equation}
It is easy to solve this system step by step. Indeed, consider 
for simplicity mono-disperse initial conditions, viz. $c_0(m)=\1_{m=1}$. Then for $t<T_{\rm gel}$ and $m=1$,  one finds $\dd c_t(1) = -c_t(1) \dd t$ and thus
$c_t(1)= \ee^{-t}$.
 Next for $m=2$ one gets by substitution
$$\dd  c_t(2) = \left(\frac{1}{2} \ee^{-2t} - 2c_t(2)\right) \dd t$$
and then $c_t(2)= \frac{1}{2}  t \ee^{-2t}$. 
By induction,  one arrives at the solution due to McLeod \cite{ML}:
\begin{equation}\label{ML}
c_t(m)= t^{m-1} m^{m-2} \ee^{-mt}/m!\,, \qquad m\in \N^*.
\end{equation}
We stress that this is only valid before the gelation time $T_{\rm gel}$,
which is a finite quantity. Indeed   $ \langle c_t, {\rm Id}\rangle$
is then given by
\begin{equation}\label{gelation}
\theta(t):=\sum_{m=1}^{\infty}\frac{t^{m-1} m^{m-1}}{m!} \ee^{-tm}\,,\end{equation}
and one can check  that $\theta(t)=1$ for $t\leq 1$ while\footnote{More precisely, observe that
 $$t\theta(t) = \sum_{m=1}^{\infty}\frac{m^{m-1}}{m!} \exp(-(t-\ln t)m)\,,$$
 so that if $t'<1<t$ are such that $t-\ln t = t'-\ln t'$, then 
$t\theta(t)=t'\theta(t')= t'.$
This entails that $\theta(t)=t'/t$ is the unique solution in $(0,1)$ to the equation $\ee^{t(\theta(t)-1)}=\theta(t)$.}
 $0<\theta(t)<1$  for $t>1$.
In particular, $T_{\rm gel}=1$.

Considering generating functions provides an alternative way for solving Equations \eqref{Smolsimple} for  $t< T_{\rm gel}$; see e.g. \cite{DeaconuTanre} and the references therein. 
Indeed, we can interpret \eqref{Smolsimple} as a convolution equation
$$\frac{1}{m}\partial_tp_t(m)=\frac{1}{2}p_t*p_t(m)-p_t(m)$$
where $p_t(m)=mc_t(m)$ is the mass-function of some probability distribution on $\N^*$. 
In terms of the generating function $g_t(x)=\sum_{m=1}^{\infty} x^m p_t(m)$, $x\in[0,1]$, this yields 
$${\partial_t} g_t(x) = \frac{x}{2}\partial_x g^2_t(x)- x\partial_x g_t(x)
= x(g_t(x)-1)\partial_x g_t(x)\,.$$
This is a quasi-linear PDE which can  be solved by the method of characteristics. One finds that for $0<t<1$ (assuming again mono-disperse initial conditions for simplicity), $g_t: [0,1]\to [0,1]$ is the inverse  bijection of  $x\mapsto x\exp(t(1-x))$. It is known that this is the generating function of the Borel distribution 
with parameter $t$, so inverting generating functions gives
$$p_t(m)=\frac{\ee^{-tm}(tm)^{m-1}}{m!}\,, \qquad m\in\N^*\,,$$
and one recovers the solution \eqref{ML} of McLeod.

We conclude this section by discussing briefly solutions after the gelation time. 
In the mono-disperse case $c_0(m)=\1_{m=1}$, 
Kokholm \cite{Kok} established global existence and uniqueness of the solution, which is given for $t\geq 1$ by 
$$c_t(m)= \frac{m^{m-2} \ee^{-m}}{m!t}\,, \qquad m\in \N^*.$$
Observe that this {\it does not} solve the simpler equation \eqref{Smolsimple} and that for $t>1$, $\langle c_t,{\rm Id}\rangle = 1/t>\theta(t)$. 
For poly-disperse initial conditions, the gelation time can be determined explicitly as
\begin{equation}\label{tgel0}
T_{\rm gel}=1/\langle c_0,m^2\rangle = 1/\sum_{m=1}^{\infty}m^2c_0(m)\,;
\end{equation}
note that gelation may occur immediately, namely $T_{\rm gel}=0$ whenever the second moment of $c_0$ is infinite.
Existence and uniqueness of the solution on the whole time interval $[0,\infty)$ have been established rigorously  only quite recently by Normand and Zambotti \cite{NZ}, although of course several partial results appeared previously in the literature; see the comments and references in \cite{NZ}.

\subsection{Connexion with the Erd\H{o}s-R\'enyi random graph model}
We then study coagulation through a microscope, in the sense that we now consider a microscopic stochastic model which follows related dynamics at the infinitesimal scale, and whose macroscopic evolution is close  to Smoluchowski's coagulation equation. 
More precisely, Marcus \cite{Marcus} and Lushnikov \cite{Lush} have introduced stochastic models  for the evolution  of finite systems of particles with total mass $n$,  in which any pair of particles, say $(m,m')$ coagulate at a rate $\kappa(m,m')/n$, independently of the other pairs. These are known as {\it Marcus-Lushnikov coalescents}; see the survey by Aldous \cite{Aldous} for background and further references. Roughly speaking, the connexion between the deterministic model of Smoluchowski and the stochastic one of Marcus and Lushnikov is  that solutions to Smoluchowski's coagulation equations can be obtained as hydrodynamic limits of Marcus-Lushnikov  processes, see  e.g.
Cepeda and Fournier \cite{CF}, Eibeck and Wagner \cite{EW},  Jeon \cite{Jeon}, Norris \cite{Norris1, Norris2}, etc. 

In the case of the multiplicative kernel $\kappa(m,m')=m\cdot m'$,  the Marcus-Lushnikov coalescent is closely related to a notorious model of random graph theory which was studied first  by Erd\H os and R\'enyi; see e.g. \cite{Bollobas} or  \cite{Remco}. Fix a number $n\gg 1$ of vertices and consider the random graph process such that  an edge connecting a given pair of vertices appears at an exponentially distributed time with parameter $1/n$ (i.e. with expectation $n$), independently of the other pairs of vertices. 
In particular, the rate at which an edge appears between a  connected component of size $m$ 
and another connected component of size $m'$ is $(m\cdot m')/n$ (this is $1/n$ times the number of edges which can possibly connect these two clusters). 
In other words, the process which records the sizes of the clusters  in the random graph process is Markovian; more precisely it is a version of the Marcus-Lushnikov multiplicative coalescent started from a mono-disperse initial condition.

Erd\H os and R\'enyi pointed at a nowadays classical phase-transition for the size of the connected component in the random graph  which contains a typical vertex. Roughly speaking, the probability that this size remains bounded when $n\to \infty$  is $\sim 1$ for $t<1$, whereas for $t>1$, this probability  is only $\sim \theta(t)<1$, where $\theta(t)$ is the unique solution in $(0,1)$ to the equation $\ee^{t(x-1)}=x$. This provides a probabilistic explanation of the gelation phenomenon; recall in particular the discussion near \eqref{gelation}. 

More precisely, the random graph model also enables us to recover the solution of McLeod before gelation. Informally, it is easily shown that for $n\gg 1$, the  cluster 
 of the random graph process at  time  $t<1$ which contains a given vertex resembles a Galton-Watson tree with reproduction law given by the Poisson distribution with parameter $t$. In particular the law of its size converges when $n\to\infty$
to the Borel distribution with parameter $t$.
Roughly speaking, this entails that  $m$ times the concentration  of clusters of size  $m$ 
(because 
a given cluster of size $m$ appears exactly $m$ times as the cluster containing a certain vertex), is close to the mass function of the Borel$(t)$ distribution evaluated at $m$, i.e. $t^{m-1} m^{m-1} \ee^{-mt}/m!$. This is the formula \eqref{ML} obtained by McLeod. 

We stress that this connection between the random graph process 
and Smoluchowski's coagulation equations for the multiplicative kernel ceases after time $1$, i.e. after the phase transition or the gelation time. The reason is intuitively clear: in the random graph process, the giant cluster which appears at time $1$ interacts thereafter with microscopic clusters
(informally, more and more microscopic clusters are absorbed by the giant component), whereas in Smoluchowski's model, the gel (i.e. giant components) does not interact with particles. 
As a matter of fact, the random graph process is even more closely related to a variation of Smoluchowski's coagulation equation, namely the Flory's equation, in which interaction with the gel is taken into account; see Fournier and Lauren\c{c}ot \cite{FL}.

\section{Coagulation equations with limited aggregations}
\subsection{Setting}
We now describe the modified model with limited aggregations; it may be useful for the intuition to think of the formation of polymers as the result of the creation of  bonds  connecting different atoms. The name {\it particle} now refers to a polymer which consist of  atoms (possibly a single one) connected by bonds. Each atom has a certain number of 
potential connexions, which we call arms,  and which can be paired with arms of other atoms to create  bonds. 
In this setting, the size $m$ of a particle is given by the number of atoms that constitute this polymer,
and its number of arms is the total number of potential connections which are still available for those atoms. 
 
 We write $c_t(a,m)$ for the concentration at time $t$ of particles with size $m$ and $a$ arms in the medium, and assume for simplicity that at the initial time, there are no bonds connecting atoms, i.e.
 $c_t(a,m)=0$ for all $m\geq 2$. We then define 
 $$\mu(a)=c_0(a,1)\,,\qquad a\in\N^*\,,$$ 
 and assume that the  first two moments are finite, i.e. 
\begin{equation}\label{norm}
A_j:=\sum_{a=1}^{\infty} a^j c_0(a,1)= \sum_{a=1}^{\infty} a^j\mu(a) <\infty\quad \hbox{for } j=1,2
 \,.
\end{equation}

 We imagine that when  two particles merge, one arm in each particle is used to perform a bond, an event which can be schematically represented as
$$(a,m) \hbox{ \& } (a',m') \ \mapsto \ ( a+a'-2, m+m')\,.$$
We further assume that the rate at which such aggregation occurs is $a\cdot a'$, i.e. it coincides with the numbers of possible connexions between these two particles. This means that the equations governing the evolution of concentrations are now
\begin{align}\label{EQ8}
\partial_t c_t(a,m) =\ &\frac{1}{2}\sum_{a'=1}^{a+1}\sum_{m'=1}^{m-1} a'
c_t(a',m')\cdot (a-a'+2) c_t(a-a'+2,m-m')&  \nonumber  \\
&-\sum_{a'=1}^{\infty}\sum_{m'=1}^{\infty}ac_t(a,m)\cdot a'c_t(a',m'), &
\end{align}
where here again the first term in the right-hand side accounts for the creation of a particle $(a,m)$ as the result of the aggregation of a particle $(a',m')$ with a particle $(a-a'+2, m-m')$, and the second for the disappearance of a particle $(a,m)$  after an aggregation event. 

Particles with no arms have a special role. They only arise from the aggregation of pairs of particles each with a single arm, and cannot be involved into  further coagulation events. The concentration $c_t(0,m)$ of such particles thus increases with time and possesses  a limit as $t\to\infty$.

\subsection{Solution and gelation}
Equation \eqref{EQ8} somewhat  resembles \eqref{Sm} for the multiplicative kernel, 
and even though it looks more involved, it can also be solved explicitly to a large extent, see \cite{Be} and \cite{NZ} for details.

Recall from \eqref{tgel0}  that the gelation time in  Smoluchowski's coagulation equations with the multiplicative kernel can be expressed as the inverse of the second moment 
of the initial condition; in particular gelation always occurs. For coagulation with limited aggregations, the gelation time is given by
$$T=  \left\{ \begin{matrix}
\infty\ \hbox{ if }A_2\leq 2 A_1, &\\ 1/(A_2-2A_1)\ \hbox{ if }A_2> 2 A_1, &\\
\end{matrix}\right.
$$
where $A_1$ and $A_2$ are defined in  \eqref{norm}. 
We stress that $T$ can be infinite, a situation which does not occur when aggregations are unlimited. The main result of \cite{Be} 
solves \eqref{EQ8} before gelation; we refer to \cite{NZ} for more complete results including post-gelation solutions. The probability measure $\nu$
on $\N$ defined by 
$$\nu(j)=\frac{j+1}{ A_1}\,\mu(j+1)\ ,\qquad j\in\N$$
plays a crucial role.

\begin{theorem}\label{T1} 
The system \eqref{EQ8} has a unique solution 
on $[0,T)$  started from
$$c_0(a,m)={\bf 1}_{\{m=1\}}\mu(a)\,,\qquad a,m\in\N^*$$
which is given for $a,m\geq 1$ by
$$c_t(a,m)= \frac{(a+m-2)!}{a! m!} A_1^m  t^{m-1}(1+A_1 t)^{-(a+m-1)}\nu^{*m}(a+m-2)\,,$$
where $\nu^{*m}$ denotes the $m$-th convolution product of the probability measure $\nu$. Further
$$c_t(0,m)=\frac{A_1}{m(m-1)}(1+1/(A_1 t))^{1-m}\nu^{*m}(m-2)\qquad \hbox{ for }m\geq 2\,.$$
\end{theorem}

Let us now  briefly sketch the main lines of the proof of Theorem \ref{T1}. Without loss of generality, we may suppose that concentrations have been normalized so that $A_1=1$. 
One starts by showing that for $t<T$,  the mean number of arms at time $t$ is simply given by
$$\sum_{a=1}^{\infty}\sum_{m=1}^{\infty} a c_t(a,m)= \frac{1}{1+t}\,,$$
and then \eqref{EQ8} can be re-written as 
$$\partial_t c_t(a,m) =\frac{1}{2}\sum_{a'=1}^{a+1}\sum_{m'=1}^{m-1} a'
c_t(a',m') (a-a'+2) c_t(a-a'+2,m-m')
-\frac{ac_t(a,m)}{1+t}. 
$$
Next one uses similar techniques of generating functions as in Section 2. That is one defines for $x,y\in[0,1]$
$$g_t(x,y)=\sum_{a=1}^{\infty}\sum_{m=1}^{\infty} x^a y^m c_t(a,m)$$
and observes that \eqref{EQ8} then implies that the function $k_t(x,y):=\partial_x g_t(x,y)$ solves for $t<T$ the quasi-linear PDE
$$\partial_t k_t(x,y)=\left(k_t(x,y)-\frac{x}{1+t}\right)\partial_x k_t(x,y)-\frac{1}{1+t}k_t(x,y)\,.$$
The latter can be solved by the method of characteristics (see e.g.  \cite{Evans} for background), and one gets
$$k_t(x,y) = (1+t)^{-1}k_0(\ell_t(x,y),y)= t^{-1} \ell_t(x,y)-\frac{x}{t^2+t}\,,$$
where $\ell_t(\cdot,y):[0,1]\to[0,1]$ is the inverse to the function  $x\mapsto (1+t)x-tk_0(x,y)$. 
All  what is needed now is to invert the generating function, which is realized thanks to Lagrange inversion formula (see, for instance, Section 5.1 in \cite{Wilf}), and we arrive at the solution of the statement.

It is interesting to point out that similar arguments also work in a more sophisticated model in which arms have genders (male or female) and only couple of arms with different genders can be paired;
see \cite{Normand}. 
\subsection{Limiting concentrations} 

 One remarkable feature of the model with limited aggregations is that concentrations have a non-trivial limit as time tends to $\infty$. In the case without gelation (i.e. $T=\infty$, or equivalently $A_2\leq 2A_1$), this is immediately seen from Theorem \ref{T1}. The case where gelation occurs is significantly harder, but somewhat surprisingly, one obtains a very similar formula. The following statement rephrases Corollary 6.2 in Normand and Zambotti \cite{NZ} (and also corrects a little misprint there);  recall that $\nu^{*m}$ has been defined in Theorem \ref{T1}.

\begin{theorem} \label{T2}  The concentration $c_t(a,m)$ of particles $(a,m)$ at time $t$ with $
a\in\N$ and $m\geq 2$, has a limit $c_{\infty}(a,m)$
as $t\to \infty$. 

\noindent{\rm (i)} When there is no gelation (i.e. $T=\infty$), the latter is given  by
$$c_{\infty}(a,m)= {\bf 1}_{\{a=0\}}\frac{A_1}{m(m-1)}\nu^{*m}(m-2)\,.$$

\noindent{\rm (ii)} When gelation occurs (i.e. $T<\infty$), let  $\eta$ denote the unique solution to
the equation $\eta g'_{\nu}(\eta)=g_{\nu}(\eta)$ where $g_{\nu}(x)=\sum x^i \nu(i)$ is the generating function of $\nu$. 
Then one has
$$c_{\infty}(a,m)= {\bf 1}_{\{a=0\}}\frac{A_1}{m(m-1)}\beta^{m-1}\nu^{*m}(m-2)\,,$$
with $\beta= \eta/g_{\nu}(\eta)=1/g'_{\nu}(\eta)>1$.
\end{theorem}
It is interesting to mention further that when gelation does not occur, we have
$$\sum_{m=2}^{\infty} mc_{\infty}(0,m)=\sum_{a=1}^{\infty}c_0(a,1)\,,$$
which means that the total mass of polymers in the terminal state coincides with the initial total mass of atoms (which is stronger than absence of gelation). 

The explicit formula for the limiting concentration has a strong probabilistic flavor. More precisely, consider a Galton-Watson process with reproduction law $\nu$, i.e.  a population model in which individuals reproduce independently one from  the other, and with the same offspring distribution $\nu$. A classical result due to Dwass \cite{Dwass} claims that the law of the size of the total population in such a Galton-Watson process with {\it two} ancestors is given by
$$\frac{2}{m}\nu^{*m}(m-2)\,,\qquad m\geq 2\,.$$
The purpose of the next section is to provide a probabilistic explanation of the appearance of Galton-Watson processes with two ancestors in this setting. 

\section{Stochastic coalescence with limited aggregations}
\subsection{The random configuration model and giant components}
The deep link between Smoluchowski's coagulation equations and Marcus-Lushnikov stochastic coalescent invites us to consider the following stochastic evolution for  finite systems of  polymers. Recall that a  particle (polymer) is characterized by its size $m$ and its number $a$ of arms; a pair of particles, say $(a,m)$ and $(a',m')$, merge at rate $a\cdot a'$ to produce a single particle 
$(a+a'-2, m+m')$, independently of the other pairs. Assume henceforth  that the initial state of the system is purely atomic, i.e. all the particles have unit size, and recall that we are mainly interested in the terminal state of the system (i.e. in the limiting concentrations for the macroscopic model). 
We  stress that arms of a polymer are only used to create bonds with other polymers;
two arms in the same polymer do not joint to form a bond. In other words,  such polymers have a tree structure.

The evolution of the system of polymers described above  is closely related to 
the random configuration model. The latter is  a simple stochastic algorithm introduced independently by Bollob\'as \cite{Bollo} and Wormald \cite{Worm} in the late 70's,  which aims at producing a random graph on a set of vertices with pre-described degrees; see Chapter 7 in \cite{Remco} for a detailed account. Typically, we start with some finite set $V$ of vertices and assign to each vertex $v\in V$ a degree $d(v)\in\N^*$. We imagine that there are $d(v)$ stubs attached to the vertex $v$, and join uniformly at random pairs of stubs to create edges between those vertices. 
The resulting multi-graph is known as the {\it random configuration model}; in general it  is not simple, in the sense that there may exist loops and multiple edges. However, it is known that, loosely speaking, if the degrees $d(v)$ are sufficiently small, then only few clusters in a random configuration have loops, multiple edges, or even cycles. If we discard this pathology, the random algorithm which produces the random configuration is then remarkably close to the stochastic coalescence with limited aggregations.  

The purpose of this section is to show that, at least in appropriate regimes for which, informally, the number of arms is not too large, results on random configurations provide a probabilistic explanation of  properties of the solutions to the deterministic coagulation equations with limited aggregations.
In this direction, we first recall that
Molloy and Reed \cite{MR1, MR2} have pointed at a phase-transition for large sparse random configurations which generalizes that of Erd\H{o}s and R\'enyi for large random graphs. Specifically,   consider for every integer $n\geq 1$ a set $V_n$ of $n$ vertices and a function
 $d_n:  V_n\to \N^*$ that specifies the number of stubs appended to each vertex.
 Introduce the empirical distribution of the number of stubs
 $$\mu_n(i):=\frac{1}{n}\#\{v\in V_n: d_n(v)=i\}\,,\qquad i\in\N^*\,.$$
 Assume  that the limit
 \begin{equation} \label{EQ7}
 \lim_{n\to\infty} \mu_n(i):=\mu(i)
 \end{equation}
 exists for every $i\geq 1$, and further  that 
 \begin{equation}\label{EQ7'}
  \lim_{n\to\infty} \sum_{i=1}^{\infty} i\mu_n(i)=\sum_{i=1}^{\infty} i\mu(i) <\infty\,.
\end{equation}
 Informally, Molloy and Reed have shown that if 
 \begin{equation}\label{Molloy}
 \sum_{i=1}^{\infty} i(i-2)\mu(i)\leq 0\,,
 \end{equation}
 then with high probability when $n\to \infty$, the size of a typical  connected component  in a random configuration on $V_n$ remains bounded, whereas when \eqref{Molloy} fails, then with high probability there exists a giant cluster of size of order $n$. 
 Observe that  \eqref{Molloy} can be rewritten as $A_2\leq 2A_1$ in the notation \eqref{norm}; so we recover the necessary and sufficient condition for the absence of gelation. 
Deeper connexions will be pointed at in the next section. 
 
 \subsection{The structure of typical polymers in the sub-critical case}
 In this section, we assume that \eqref{EQ7}, \eqref{EQ7'} and  \eqref{Molloy} holds, so typical clusters in the random configuration model remain microscopic as the number of vertices tends to infinity, and gelation does not occur. Our purpose is to show that one can recover the formula of Theorem \ref{T2} for the limiting concentration, by statistical analysis of those microscopic clusters. 
This was achieved in \cite{BS}, here we shall merely provide a rather informal sketch of the argument, refereeing to that article for details.

Consider the random configuration model on a set $V_n$ of $n$ vertices. For every {\it oriented} edge $e$ of this random graph, we consider the connected component $C_e$ that contains this edge, viewed as a combinatorial structure on a set of unlabeled  vertices (in the sense that two isomorphic subgraphs are identified) with a distinguished oriented edge. We are interested in the empirical measure
 $$\epsilon_n(\Gamma)= \frac{1}{S_n}\#\{e: C_e=\Gamma\}$$
 where $S_n$ is the total number of arms (or equivalently $S_n$ is the number of oriented edges) and $\Gamma$ is a generic combinatorial structure with a distinguished oriented edge. 
 
Recall  that $\nu$ denotes the probability measure given by
$$\nu(i) = (i+1)\mu(i+1)/A_1\,,\qquad i\in \N\,.$$
The condition  \eqref{Molloy} is equivalent to requiring that $\nu$ is critical or sub-critical in the sense of branching processes, that is $\sum i \nu(i)\leq 1$. 
Now consider a Galton-Watson process started from two ancestors, with reproduction law $\nu$. We know that the process becomes extinct a.s., and we can consider the two genealogical trees of these two ancestors, which we further connect by an oriented edge from the first ancestor to the second one. This gives a combinatorial structure with one distinguished oriented edge; we write 
$\mathbf{GW}_{\nu,2}$ for its law.

The main result of \cite{BS} is the following:

\begin{theorem}\label{T3} Suppose \eqref{Molloy} holds. Then 
$$\lim_{n\to\infty}\epsilon_n = \mathbf{GW}_{\nu,2}\qquad\hbox{in probability,}$$
where the space of probability measures on the discrete set of combinatorial structures with one distinguished edge is endowed with the total variation distance. 
\end{theorem}

Let us briefly present the idea of the proof, which already appears in different forms in the literature on random configurations. From the construction of the configuration model, we see that sampling an oriented edge uniformly at random amounts to picking two arms uniformly at random at the initial stage of the construction (i.e. before any bond between vertices has been formed) and joining them. Let $v$ and $v'$ be the vertices to which these two arms are attached; possibly $v=v'$, but this event is quite unlikely when $n\gg 1$ and will thus be discarded. Note the degree of $v$, $d_n(v)$
does not have the law  $\mu_n$ normalized to be a probability measure, but rather the size-biased version of the latter.  It follows from our assumptions that
$d_n(v)-1$, which represents the number of arms of $v$ which remain available after $v$ and $v'$ get connected, converges in distribution to $\nu$, and by an obvious symmetry, $d_n(v')-1$ also converges in distribution to $\nu$. More precisely, $d_n(v)-1$ and $d_n(v')-1$ are asymptotically independent, in the sense that the pair $(d_n(v)-1,d_n(v')-1)$ converges in law   to a pair of independent variables with the same law $\nu$.

Then let us work conditionally on  $d_n(v)-1=k$ and $d_n(v')-1=\ell$. 
Informally, each available arm of $v$ and of $v'$ will then be connected to different
vertices, say $v_1, \ldots, v_k$ and $v'_1, \ldots, v'_{\ell}$ such that, when $n\to \infty$, the variables
$d_n(v_1)-1, \cdots, d_n(v_k)-1, d_n(v'_1)-1, \cdots, d_n(v'_{\ell})-1$ become asymptotically independent and have all the same distribution $\nu$. By iteration, we see that if we pick an oriented edge $e$ uniformly at random, the cluster $C_e$ converges in distribution to a random combinatorial structure with law $\mathbf{GW}_{\nu,2}$. 

The proof is then completed with an argument of propagation of chaos. Namely, we have to check that if we now pick two oriented edges uniformly at random, say $e$ and $e'$, then the 
combinatorial structures $C_e$ and $C_{e'}$ become asymptotically independent when $n\to \infty$.  Indeed, it is well-known that propagation of chaos is equivalent to the weak convergence of the empirical distributions; see \cite{ASz}.

We have now all the ingredients needed for a probabilistic (or microscopic) explanation of Theorem \ref{T2}(i). First, the total number of arms in the system is
$$S_n= n\sum_{1}^{\infty} i \mu_n(i) \sim n \sum_{1}^{\infty} i \mu(i)= nA_1.$$
Then, roughly speaking, we know from Theorem \ref{T3} that most clusters in the random configuration model are trees; 
in particular a tree-cluster $C$ of size $|C|=m$ (number of vertices) has exactly $m-1$ edges, which in turn correspond to $2(m-1)$ arms. If we  think of the random configuration as the terminal state of the polymer system, then
the concentration of polymers of size $m$ is  approximatively
$$c_{\infty}(m)\approx \frac{S_n}{2(m-1)n} \epsilon_n(\{C: |C|=m\})$$
Recall also  that thanks to Dwass' formula,
$$\mathbf{GW}_{\nu,2}(|C|=m)= \frac{2}{m} \nu^{*m}(m-2)\,,$$
so combining with Theorem \ref{T3}, we arrive at
$$c_{\infty}(m)\approx \frac{A_1}{m(m-1)} \nu^{*m}(m-2)\,,$$
 which is the formula found in Theorem \ref{T2}. 

\subsection{Formation of gel and self-organized criticality}
We now conclude this text by presenting succinctly and informally some resent results obtained by Merle and Normand (see Chapter 4 in \cite{MN}) in the supercritical case. 

Roughly speaking,  stochastic coalescence can be related to Smoluchowski's coagulation equations only up to the time when a giant component appears, which corresponds to the gelation time. This stems from the fact that in the stochastic framework, giant components interact with smaller components, whereas this is not the case in  Smoluchowski's setting.
Fournier and Lauren\c{c}ot \cite{FL} pointed at a modified version of stochastic coalescence which 
circumvents this problem. Denoting by $M$  the total mass of particles in the stochastic system, they introduce a threshold $\alpha(M)$ with $1\ll \alpha(M)\ll M$,  such that as soon as a particle
with size greater than $\alpha(M)$ appears, it falls into the gel and ceases to interact with smaller particles. In other words,  particles with size greater than $\alpha(M)$ are seen as giant, 
and we may imagine that they are then removed from the system.  Thus only particles with mass less than $\alpha(M)$ are allowed to coagulate, and an hydrodynamic limit procedure now enables one to relate  this modified stochastic coalescence to Smoluchowski's coagulation equations for all times. Note also this modified model coincides with the original one of Marcus and Lushnikov before the appearance of giant particles.

Merle and Normand exploited this idea for stochastic coalescence with limited aggregations.
Strictly speaking, they consider a slightly different model in which time runs at a different speed, however this is only a minor change made for convenience; in particular this does not alter the terminal state of the system. So imagine that we start with $n$ atoms with arms. As time passes, polymers are formed by the creation of bonds connecting atoms just  as in Section 4.1, except that now polymers with size greater than the threshold\footnote{Merle and Normand assume that $n^{1/3}\ll \alpha(n) \ll n$, which is slightly less general than the condition of Fournier and Lauren\c{c}ot.}
$\alpha(n)$ fall instantaneously into the gel.

Merle and Normand focus on the total concentration, say $m_t^n$, of atoms (not polymers) in solution at time $t$, i.e. which have not fallen into the gel before time $t$,  and on the empirical measure $\pi_t^n$ of used arms in the solution. That is $n m_t^n$ represents the total number of atoms in solution at time $t$, and the number of atoms in solution at time $t$ for which exactly $k$ arms have already been used is $n m_t^n \times \pi_t^n(k)$. They point at the following simple description of the distribution of the configuration in solution: conditionally on $n_t^n$ and $\pi_t^n$, the configuration in solution has the same law as a random configuration model started with $n_t^n$ atoms and empirical distribution of arms $\pi_t^n$, and conditioned on having no cluster of size greater than $\alpha(n)$.

They prove that under hypotheses a bit more restrictive than 
\eqref{EQ7} and \eqref{EQ7'}, $m_t^n$ converges as $n\to \infty$ to some deterministic $m_t\in[0,1]$ and $\pi_t^n$ converges to some deterministic probability measure on $\N^*$. Both $m_t$ and $\pi_t$ can be characterized explicitly in terms of the limiting distribution of arms $\mu$. It should not come as a suprise that before the gelation time, $\pi_t$ is subcritical in the sense of the random configuration model, i.e.
$$\sum_{i=1}^{\infty} i(i-2)\pi_t(i)<0$$
(recall the condition \eqref{Molloy} of Molloy and Reed); but a much more interesting feature is that
$\pi_t$ remains {\it exactly} critical after the gelation, i.e. 
$$\sum_{i=1}^{\infty} i(i-2)\pi_t(i)=0\,.$$
This provides a nice illustration of {\it self-organized criticality}, a general phenomenon which is often observed in Statistical Physics (see \cite{BTW}), but has been rarely be proven mathematically.

Finally, Merle and Normand consider the limits $m_{\infty}=\lim_{t\to \infty} m_t$
and $\pi_{\infty}=\lim_{t\to\infty} \pi_t$. Specifically, they obtain that
$$m_{\infty}=\sum_{i=1}^{\infty} \eta^i \mu(i)\,,$$
where $\eta$ has been defined in Theorem \ref{T2}(ii), 
and that 
$$\pi_{\infty}(i)= \frac{\eta^i} {m_{\infty}} \ \mu(i)\,,\qquad i\in\N^*\,.$$
This enables them to repeat the argument at the end of Section 4.2 and to provide a probabilistic explanation for the simple formula of Theorem \ref{T2}(ii) which gives the limiting concentrations in the supercritical case.

\end{document}